\documentclass[a4paper,11pt]{article}
\usepackage[retainorgcmds]{IEEEtrantools}
\usepackage{algorithmic,algorithm}
\usepackage{epsfig}
\usepackage{array}
\usepackage{amsmath}
\usepackage{amssymb}
\usepackage{amsthm}
\usepackage{color}

\usepackage{cite}

\usepackage{xcolor,calc}

\title{Policy-based reserves for power systems}
\author{Joseph Warrington, Paul Goulart, S\'{e}bastien Mari\'{e}thoz, Manfred Morari%
\thanks{This work was supported by the European Commission Research Project FP7-ICT-249096 \emph{Price-based Control of Electrical Systems (E-PRICE)}.}
\thanks{The authors are with the Automatic Control Laboratory, Swiss Federal Institute of Technology (ETH) Zurich,
Physikstrasse 3, 8092 Zurich, Switzerland. Contact: {\tt\small warrington@control.ee.ethz.ch}}}%

\begin{document}
\maketitle

\begin{abstract}
This paper introduces the concept of affine reserve policies for accommodating large, fluctuating renewable infeeds in power systems. The approach uses robust optimization with recourse to determine operating rules for power system entities such as generators and storage units. These rules, or policies, establish several hours in advance how these entities are to respond to errors in the prediction of loads and renewable infeeds once their values are discovered. \emph{Affine} policies consist of a nominal power schedule plus a series of planned linear modifications that depend on the prediction errors that will become known at future times. We describe how to choose optimal affine policies that respect the power network constraints, namely matching supply and demand, respecting transmission line ratings, and the local operating limits of power system entities, for all realizations of the prediction errors. Crucially, these policies are time-coupled, exploiting the spatial and temporal correlation of these prediction errors. Affine policies are compared with existing reserve operation under standard modelling assumptions, and operating cost reductions are reported for a multi-day benchmark study featuring a poorly-predicted wind infeed. Efficient prices for such ``policy-based reserves'' are derived, and we propose new reserve products that could be traded on electricity markets.
\end{abstract}

\section{Introduction}

A key challenge in incorporating highly variable intermittent renewable energy sources into power systems is the need to maintain system integrity while making best use of the energy they provide, which comes at zero marginal cost. It is widely agreed that in the next few decades, as the share of wind power becomes very large, current techniques for accommodating wind variability will become sufficiently expensive that alternatives will be sought \cite{ackermann_european_2007, newbery_market_2010}.

Running power systems with very high wind penetration and without excessive frequency control costs, or resorting to curtailment of renewable output, requires intelligent use of the best available forecasts, at all times. In particular, any successful method for dealing with the high variability of renewables on intraday timescales (the only timescales over which prediction errors are reasonably small \cite{giebel_forecast_2007}) will require the following: 
\begin{enumerate}
\item A forecast of future intermittent energy injections available at the time when control decisions are made.
\item Rules for acting on errors in this forecast when they are discovered.
\item Forecast error probability distributions and their correlations in both time and space over the grid.
\end{enumerate}

Theoretical attention to these points has grown in the last few years as the share of wind power in several countries has grown \cite{galiana_scheduling_2005, bouffard_stochastic_2008, morales_economic_2009, xiao_operating_2011, vrakopoulou_probabilistic_2012}. In Morales et al.~\cite{morales_economic_2009} a unit commitment integer programming problem was solved first, and then in a second stage reserve margins were selected based on the requirement that the actual reserve deployment be feasible for all the scenarios considered. These ideas have since been developed to provide probabilistic guarantees on transmission constraint satisfaction using a limited number of scenarios \cite{vrakopoulou_probabilistic_2012}.

In this paper, we consider how reserves could be operated more efficiently on a daily timescale around any unit commitment decisions that have already been made. In contrast with existing literature, we use a robustness formulation of the problem where the bounds on the uncertainty are assumed to have been chosen according to probabilistic criteria or are inherent, e.g. arising from wind farm capacities. We choose optimal time-coupled policies, which are rules agreed in advance governing how individual power system entities will respond to prediction errors affecting power system operation. Policies therefore constitute a reserve mechanism, a concept that represents the main contribution of this work. We report reductions in the cost of operating reserves when policies are used in the presence of a large uncertainty, in comparison to reserve rules that do not make use of policies.
 

This work was inspired by results on disturbance feedback policies from the control literature. However these are predated by the concept of linear decision rules (LDRs) from operations research, where current states, past data or future predictions are combined linearly in order to make an operational decision \cite{holt_linear_1955,silver_simple_1979}. Typically, though, LDRs were unable to deal rigorously with operating constraints, and were not studied in much detail after the 1970s \cite{kuhn_primal_2009}. 

In the last decade, however, LDRs have been revived as a means of solving \emph{constrained} optimization problems where the minimizer is allowed to be a function of the data uncertainty \cite{guslitser_uncertainty-immunized_2002, ben-tal_adjustable_2004, kuhn_primal_2009}. This has led to some new applications, for example in portfolio optimization \cite{rocha_multistage_2012}. These solution methods were also shown to be applicable to robust predictive control \cite{lofberg_minimax_2003, goulart_optimization_2006}, where control policies with various dependences on the disturbance have been studied as a means of respecting state and input constraints under uncertain system dynamics. A recent application of this is intelligent building control \cite{oldewurtel_use_2012}.

Although optimal short-term operation of power systems, including reserves, has been studied in various ways for decades \cite{bertsekas_optimal_1983,wang_short-term_1995} as a variant of the standard optimal power flow problem, affine policies have not until now been exploited for real-time decision making in electricity provision under uncertainty, despite their attractiveness for incorporating forecasts into power system operations. To this end, we present systematic ways of using future wind prediction error measurements to reduce the average costs of reserve provision. We use a linearized transmission model, which has been consistently shown in many real applications to produce good approximations of true AC power flows.



\subsection*{Summary of notation by theme}
~\\
\begin{tabular}{ll}
$(\cdot)'$ & Vector or matrix transpose \hfill \\
$[\cdot]_k$ & Element $k$ of a vector \\
$\langle X,Y \rangle$ & Trace of product $X'Y$ \\
$\otimes$ & Kronecker product \\
$I_n$ & Identity matrix of dimension $n$
\medskip	 \\
$T$ & Length of time horizon in steps \\	
$N_{\rm p}$ & Number of participants \\
$r_i$ & Nominal inelastic power flow from participant $i$ \\
$G_i$ & Map from uncertainty to inelastic power flows \\
$\delta_k$ & Uncertain vector at time $k$ \\
$\delta$ & Stacked vector of future uncertain quantities \\
$N_\delta$ & Elements of uncertainty vector per time step \\
$\Delta$ & Set from which uncertainty is drawn \\
$S,h$ & Parameters defining $\Delta$ via inequalities \\
$q$ & Number of rows in $S$ and $h$
\medskip \\
$x^i_k$ & State of participant $i$ at time $k$ \\
$u^i_k$ & Input to participant $i$ at time $k$ \\
$n_i$   & Dimension of participant $i$'s state \\
$\tilde{A}_i$ & State transition matrix for participant $i$ \\
$\tilde{B}_i$ & Input transition matrix for participant $i$ \\
\end{tabular} 
~\\
\begin{tabular}{ll}
$\mathbf{x}^i$ & Stacked vector of participant $i$'s future states \\
$\mathbf{u}^i$ & Stacked vector of participant $i$'s future inputs \\
$A_i$ & Stacked state transition matrix for participant $i$ \\
$B_i$ & Stacked state transition matrix for participant $i$ \\
$C_i$ & Stacked output matrix for participant $i$
\medskip \\
$J_i(\cdot,\cdot)$ & Cost function for participant $i$ \\
$f_i^x$; $H_i^x$ & Linear; quadratic state cost coefficient \\
$f_i^u$; $H_i^u$ & Linear; quadratic input cost coefficient \\
$c_i$ & Constant cost component \\
$\mathcal{Z}_i$ & Local constraint set for participant $i$ \\
$T_i,U_i,V_i,w_i$ \hspace{-0.9cm} & \hspace{0.5cm} Parameters defining $\mathcal{Z}_i$ \\
$l_i$ & Number of linear inequalities defining set $\mathcal{Z}_i$
\medskip \\
$N_{\rm n},L $ & Number of transmission network nodes, lines \\
$\Gamma_i$ & Line flow contribution factor of participant $i$ \\
$\overline{p}$ & Stacked vector of line flow constraints
\medskip \\
$\pi_i(\delta)$ & General control policy for participant $i$ \\
$D_i$ & Matrix adjusting power flows in response to $\delta$ \\ 
$e_i$ & Participant $i$'s nominal elastic power flow \\
$\mathcal{F}_i(x_0^i)$ & Set of feasible policies at current state $x_0^i$ \\
$\tilde{J}_i(\cdot,\cdot,\cdot)$ \hspace{-0.6cm} & Expected cost for given reserve policy 
\medskip \\
$Z$ & Auxiliary matrix for reformulation of \eqref{eq:AllPolLineConstraint} \\
$Y_i$ & Auxiliary matrix for represention of set $\mathcal{Z}_i$ \\
$\lambda,\Pi,\nu,\Psi$\hspace{-0.5cm} & \hspace{0.2cm} Lagrange multipliers for \eqref{eq:DistrConstr1} to \eqref{eq:DistrConstr4} \\
$\lambda_i$ & Vector of power prices seen by participant $i$ \\
$\Pi_i$ & Reserve policy prices seen by participant $i$ 
\medskip \\
$q_k$ & State of random process driving uncertainty \\
$\beta_k$ & Random transition in random process \\
$\Sigma; A_\beta,b_\beta$ \hspace{-0.5cm} & \hspace{0.2cm} Variance; bounding parameters governing $\beta_k$ \\
$k$ & Time index for optimization variables \\
$t$ & Time index used for power system simulations \\
$x^i(t)$ & Realized state of participant $i$ at time $t$ \\
$u^i(t)$ & Realized input of participant $i$ at time $t$ \\
\end{tabular}


\section{Power system model} \label{sec:PSModel}
This paper considers the problem of optimal operation of an electrical network to satisfy loads in the presence of uncertainty. The uncertainty to be accommodated manifests itself in the form of random power infeeds from renewables and fluctuating load requirements. We will choose operating rules that apply only for a finite time into the future, on the assumption that new rules will be determined before the chosen rules expire. This finite time, or planning horizon, is divided into $T$ discrete time steps, corresponding to the trading intervals (of length 5 minutes to 1 hour) over which electricity is traded on modern intra-day markets \cite{stoft_power_2002}. The length of time horizon considered relevant to this work is up to 24 hours, after which predictions of renewable infeed are assumed to become too poor to incorporate into sophisticated decision-making rules, and unit commitment decisions are not yet fixed. A small worked example using the model outlined below can be found in our earlier paper \cite{warrington_robust_2012}.

\subsection{Participant model} \label{sec:PartModel}

We consider the actions of $N_{\rm p}$ generic entities, or participants, connected to a transmission grid. Each participant $i$, for example a generator, load, or storage unit, injects power into or extracts power from a fixed location on the network, in two forms (one or both of which may be present for a given participant):
\begin{itemize}
	\item \emph{Inelastic power flows}, which cannot be influenced by control signals.
	\item \emph{Elastic power flows}, which are determined by the result of an optimization over possible control actions.
\end{itemize}

\subsubsection*{Inelastic power flows}

The inelastic, or exogenous, injection or extraction of power for each participant $i$ is modelled as $r_i + G_i\delta$, with positive values denoting a net power injection. Its two components are a nominal prediction $r_i \in \mathbb{R}^{T}$ plus a linear function $G_i \in \mathbb{R}^{T \times N_\delta T}$ of entries of a random forecast error vector $\delta \in \mathbb{R}^{N_\delta T}$, whose value is to be discovered in the future. It has the form $\delta = [\delta_0',\,\ldots , \delta_{T-1}']'$, where each $\delta_k \in \mathbb{R}^{N_\delta}$. In other words, $N_\delta$ is the number of elements in the disturbance vector at a given time, and this vector is mapped to the exogenous power flows in the system at that time. If the prediction error $\delta$ turns out to be zero, then the net power injection at step $k$ will simply be $[r_i]_k$.

The forecast error $\delta$ is assumed to belong to a compact set $\Delta := \{\delta \, | \, S \delta \leq h\}$ with $h \in \mathbb{R}^q$, whose interior contains the origin. Although the error is random, the mean prediction error $\mathbb{E}[\delta]$ and the second moment $\mathbb{E}[\delta \delta'] \in \mathbb{R}^{N_\delta T \times N_\delta T}$ are assumed to be known. No other restrictions are placed on the probability distribution.

\subsubsection*{Elastic power flows}

Elastic power flows are governed by a participant's dynamics in conjunction with some pre-defined costs. We describe them in the standard state space form from systems and control (see \cite{maciejowski_predictive_2002}, Section 2.1). At time $k$, each participant $i$ has internal state $x^i_k \in \mathbb{R}^{n_i}$, where $n_i$ is the state dimension, and is governed by linear time-invariant dynamics, so that given an input $u^i_k$ at time $k$ the state at time $k+1$ is given by $x^i_{k+1} = \tilde{A}_ix^i_k + \tilde{B}_iu^i_k$, where $\tilde{A}_i \in \mathbb{R}^{n_i \times n_i}$ and $\tilde{B}_i \in \mathbb{R}^{n_i}$. The first element $[x^i_k]_1$ of the state vector $x_k^i$ is assumed to represent the current power injection at the relevant node of the transmission network, and other elements are used to model internal dynamics or memory of previous states. The \emph{scalar} input $u^i_k \in \mathbb{R}$ controls the net power injection of the participant at time $k+1$.

Assigning the current time the value $k=0$, a vector of future states for participant $i$, $\mathbf{x}^i := [x^i_1{}'\,\ldots\,x^i_T{}']' \in \mathbb{R}^{n_iT}$ can be written as a function of the input sequence $\mathbf{u}^i := [u^i_0{}\,\ldots\,u^i_{T-1}{}]' \in \mathbb{R}^{T}$ and the current state $x^i_0$:
\begin{equation}
\mathbf{x}^i = A_i x^i_0 + B_i \mathbf{u}^i \, , \label{eq:Dynamics}
\end{equation}
where
\begin{displaymath}
A_i := \left[ \begin{smallmatrix} \tilde{A}_i \\ \tilde{A}_i^2 \\ \vdots \\ \tilde{A}_i^T \end{smallmatrix} \right], \quad B_i := \left[ \begin{smallmatrix} \tilde{B}_i & 0 & \cdots & 0 \\ \tilde{A}_i\tilde{B}_i & \tilde{B}_i & \ddots & 0 \\ \vdots & \ddots & \ddots & 0 \\ \tilde{A}_i^{T-1}\tilde{B}_i & \cdots & \tilde{A}_i\tilde{B}_i & \tilde{B}_i \end{smallmatrix} \right].
\end{displaymath}
 
The vector of outputs $C_i\mathbf{x}^i$ as seen by the network is just the power injected or consumed by the participant. Therefore each matrix $C_i \in \mathbb{R}^{T \times n_i T}$ selects only the first element of the state vector at each time, i.e.\ $C_i = I_T \otimes \tilde{C}_i$, where $\tilde{C}_i = [1\,\, 0_{1 \times (n_i-1)}]$.
\subsubsection*{Costs}

The function $J_i :\mathbb{R}^{n_iT} \times \mathbb{R}^{T} \rightarrow \mathbb{R}$ is used to define costs for the states and inputs along the time horizon. We make the common assumption \cite{momoh_review_1999} that costs can be modelled using a quadratic function, 
\begin{equation}
J_i(\mathbf{x}^i,\mathbf{u}^i) := f_i^x{}'\mathbf{x}^i + \frac{1}{2}\mathbf{x}^i{}'H_i^x \mathbf{x}^i + f_i^u{}'\mathbf{u}^i + \frac{1}{2}\mathbf{u}^i{}'H_i^u \mathbf{u}^i + c_i \, ,\label{eq:CostFunction}
\end{equation}
where the Hessian matrices $H_i^x$ and $H_i^u$ are assumed to be positive semi-definite, in order for the optimization problem defined in Section \ref{sec:FeedbackPols} to be convex. Linear components $f_i^x$ and $f_i^u$ are of the form $\mathbf{1}_{T\times 1} \otimes \tilde{f}_i^x$ and $\mathbf{1}_{T\times 1} \otimes \tilde{f}_i^u$ respectively, where $\tilde{f}_i^x \in \mathbb{R}^{n_i}$ and $\tilde{f}_i^u \in \mathbb{R}$. Similarly the quadratic components are given by $H_i^x := I_T \otimes \tilde{H}_i^x$ and $H_i^u := I_T \otimes \tilde{H}_i^u$, where $\tilde{H}_i^x \in \mathbb{R}^{n_i \times n_i}$ and $\tilde{H}_i^u \in \mathbb{R}$. Coupling costs between time steps can be represented by augmenting the state vector to include a memory of prior states in the state vector. Constant $c_i := T\tilde{c}_i$ allows for a constant stage cost $\tilde{c}_i$.

\subsubsection*{Constraints}

The set $\mathcal{Z}_i$ consists of permissible combinations of state and input sequences $\mathbf{x}^i$ and $\mathbf{u}^i$ for participant $i$, which may in some cases be additionally constrained by $\delta$. It is a compact set defined by $l_i$ linear inequalities (i.e.~a polytope) and takes the form
\begin{equation}
	\mathcal{Z}_i := \left\{ \left[ \left. \begin{array}{c} \mathbf{x}^i \\ \mathbf{u}^i \\ \delta  \end{array} \right] \, \right|  \, T_i \mathbf{x}^i + U_i \mathbf{u}^i + V_i \delta  \leq w_i \right\} \, , \label{eq:ZSetDef}
\end{equation}
where $T_i \in \mathbb{R}^{l_i \times n_i T}$, $U_i \in \mathbb{R}^{l_i \times T}$, $V_i \in \mathbb{R}^{l_i \times N_\delta T}$ and $w_i \in \mathbb{R}^{l_i}$. Of course, $\mathbf{x}^i$ and $\mathbf{u}^i$ are related by the dynamic equation \eqref{eq:Dynamics}, and $T_i$, $U_i$, $V_i$, and $w_i$ may also depend on the current state $x_0^i$; these dependences are left out of the notation above for clarity. A generator with input limited to the range $[p_\text{min},p_\text{max}]$, for instance, could be modelled with $T_i=0$, $U_i =[{}^{\,\,I}_{-I}]$, $V_i = 0$, $w_i = [{}^{\,\,p_\text{max} \cdot \mathbf{1}}_{-p_\text{min} \cdot \mathbf{1}}]$. Note that since $\mathbf{x}^i$ and $\mathbf{u}^i$ are trajectories rather than state or input vectors corresponding to a single time, a wide range of constraints coupling states and inputs can be modelled. For example, ramp rates may be imposed on generators, and empty/full constraints may be imposed on storage units. 

Usually $V_i = 0$, unless the uncertainty feeds directly into the participant's operating constraints. An example of this would be a curtailable wind farm whose maximum power availability at any given time is uncertain, and whose power output could be varied at any time between zero and this upper limit. In this paper, though, wind curtailment is assumed to be undesirable.

Binary decision variables, which would be needed to model start-up and shut-down (unit commitment) decisions for generators, are for the sake of clarity not considered in this paper, since we assume that such decisions have been made at an earlier stage. However it would be possible to include integer switching inputs into the formulation we give here by using an approach such as the mixed logical-dynamical system description \cite{bemporad_control_1999}. This would lead to a mixed-integer optimization rather than the convex problem \eqref{eq:DistrOpt} arrived at in this paper. The result would be a set of unit commitment decisions taking the availability of affine policies into account. In related work, \cite{bertsimas_adaptive_2013} considers how to compute an adaptive unit commitment under uncertainty, but the formulation does not incorporate correlation information on the uncertainty, and in contrast to our method minimizes a worst-case linear cost.

\subsection{Network model} \label{sec:NetConstraints}

The network model is a standard linearized approximation of a high-voltage transmission grid \cite{christie_transmission_2000}, in which lines are lossless, voltage magnitudes are constant, and line flows are proportional to the phase differences (assumed to be small) between nodal voltages. Let each participant be connected to one of $N_{\rm n}$ network nodes, and let $L$ be the number of lines connecting these nodes.

The network imposes two constraints on power system operation. The first is that the \emph{net} power injection, comprising the sum of inelastic flows $[r_i + G_i \delta]_k$ and elastic flows $[C_i\mathbf{x}^i]_k$, has to be zero at all times $k=1,\ldots,T$. This can be modelled using an equality constraint with $T$ rows:
\begin{equation}
	\sum_{i=1}^{N_{\rm p}} (r_i + G_i \delta + C_i\mathbf{x}^i) = 0\, . \label{eq:PowerMatchConstraint}
\end{equation}

The second constraint is that line currents cannot exceed the respective line ratings anywhere on the network, at any time. Under the preceding assumptions, this constraint is linear in the power injections, as long as the net power injection into the network is zero (i.e.~condition \eqref{eq:PowerMatchConstraint} holds) \cite{christie_transmission_2000}. It can be represented by the vector inequality
\begin{equation}
	\sum_{i=1}^{N_{\rm p}} \Gamma_i(r_i + G_i \delta + C_i\mathbf{x}^i) \leq \overline{p}\, . \label{eq:LineLimitConstraint}
\end{equation} 
This has $2LT$ rows, one for each flow direction, for each line, at each time. Each matrix $\Gamma_i \in \mathbb{R}^{2LT \times T}$ maps the power output of the node to which participant $i$ is attached to contributions to line flows. Each $\Gamma_i$ can be constructed from the network line impedances using the derivation of equation (III.4) in \cite{christie_transmission_2000}, since we model line constraints as limits on the phase angle differences between adjacent transmission buses.


\section{Choosing optimal reserve policies} \label{sec:FeedbackPols}

Current reserve mechanisms use a cascaded loop structure, where the fastest (primary) controller stabilizes the grid frequency, a minutes-scale (secondary) controller corrects it back to its reference, and slower, separately-purchased tertiary reserves redispatch generators in order to free up the margins within which the faster control operates \cite{stoft_power_2002}. This reserve action is only a real-time response to the error as it unfolds, and makes no systematic use of what is expected to happen in future. Intraday electricity markets in many countries are currently experiencing dramatically increasing trade volumes \cite{weber_adequate_2010}, and this increase can be seen as an attempt to adjust the short-term economic operation of the power system in the light of new forecast information. Different countries operate these markets in diverse ways (an overview of the various mechanisms used to acquire reserves and short-term power commitments in European countries can be found in \cite{entso-e_survey_2012}). Such trading actions take little systematic account of time-coupled costs and constraints imposed on the market participants.

In this section we describe a more systematic predictive mechanism that explicitly takes account of short-term uncertainties with the aim of reducing the expected running costs of the power system over the time horizon.

\subsection{Finite horizon optimization} 

Consider the problem of minimizing expected running costs $\sum_{i=1}^{N_{\rm p}}\mathbb{E}[J_i(\mathbf{x}^i,\mathbf{u}^i)]$ over a horizon of length $T$, subject to the local constraints \eqref{eq:ZSetDef} and network constraints \eqref{eq:PowerMatchConstraint} and \eqref{eq:LineLimitConstraint}. We do this by choosing a sequence of control inputs $\mathbf{u}^i$ for each participant $i$ that can vary with $\delta$. We wish to choose the best causal response to prediction errors, a \emph{policy} $\mathbf{u}^i = \pi_i(\delta)$, where $\pi_i : \mathbb{R}^{N_\delta T} \rightarrow \mathbb{R}^{T}$ is to be chosen before the error is known. ``Causal'' means that $u^i_m$ can depend only on the measurements of $\delta_0,\delta_1,\ldots,\delta_m,\delta_{m+1}$. That is, we assume that $\delta_{m+1}$, the sub-vector of $\delta$ pertaining to time $m+1$, is revealed just \emph{before} input $u^i_m$ is applied. Obviously, a dependence on any of $\delta_{m+2},\ldots,\delta_{T-1}$ would violate causality because $u^i_m$ would be a function of information unavailable at time $m$.

Substituting $\mathbf{u}^i = \pi_i(\delta)$ into the state update equation \eqref{eq:Dynamics} and eliminating $\mathbf{x}^i$, we obtain the following finite horizon optimization problem:
\begin{subequations} \label{eq:OptOverAllPolicies}
\begin{IEEEeqnarray}{L}
\min_{\text{Causal } \pi_i} \quad  \sum_{i=1}^{N_{\rm p}}\mathbb{E}[J_i(A_ix_0^i + B_i\pi_i(\delta),\pi_i(\delta))] \label{eq:AllPolObjFunction}\\
\text{s.t.}\,  \sum_{i=1}^{N_{\rm p}} r_i + G_i \delta + C_i(A_ix^i_0 + B_i\pi_i(\delta)) = 0 \, ,  \forall \delta \in \Delta \, , \label{eq:AllPolMatchConstraint}\\
\sum_{i=1}^{N_{\rm p}} \Gamma_i (r_i + G_i \delta + C_i(A_ix^i_0 + B_i\pi_i(\delta))) \leq \overline{p} ,  \forall \delta \in \Delta,  \label{eq:AllPolLineConstraint}\\
\left[ \begin{array}{c} A_ix_0^i + B_i\pi_i(\delta) \\ \pi_i(\delta) \\ \delta  \end{array} \right] \in \mathcal{Z}_i \, , \forall \delta \in \Delta \, . \label{eq:AllPolZConstraint}
\end{IEEEeqnarray}
\end{subequations}
Constraints \eqref{eq:AllPolMatchConstraint} and \eqref{eq:AllPolLineConstraint} are the expanded forms of \eqref{eq:PowerMatchConstraint} and \eqref{eq:LineLimitConstraint} after substituting definitions of $\mathbf{x}^i$ and $\mathbf{u}^i$, so that the optimization is expressed only in terms of $\pi_i$.

This problem is intractable due to the wide variety of candidate functions $\pi_i$ that could satisfy the constraints. We therefore restrict ourselves from now on exclusively to policies of the affine form
\begin{equation}
\mathbf{u}^i = D_i\delta + e_i\, , \label{eq:AffinePolDef}
\end{equation}
so that participant $i$'s power schedule $\mathbf{u}^i$ is parameterized by a nominal schedule $e_i = [e^i_0 \ldots e^i_{T-1}]'$ plus a linear variation $D_i$ with future prediction errors. In order for the use of future disturbances to be causal, $D_i$ takes the lower-triangular form
\begin{displaymath}
D_i =  \left[ \!\begin{smallmatrix} [D_i]{}_{0,0} & 0           & \cdots          & 0 \\  
                                 \text{$[D_i]{}_{1,0}$} & [D_i]_{1,1} & \ddots          & \vdots \\ 
                                 \vdots        & \ddots      & \ddots          & 0 \\ 
                                 \text{$[D_i]{}_{T-1,0}$} & \cdots      & [D_i]{}_{T-1,T-2} & [D_i]{}_{T-1,T-1} \end{smallmatrix} \! \right]
\end{displaymath}
where $[D_i]_{l,m} \in \mathbb{R}^{1 \times N_\delta}$ is the response of input $u^i_l$ to error $\delta_{m+1}$. The presence of non-zero elements above the diagonal would violate causality, because the values of as-yet-unknown errors would contribute to the control rule.

Constraint \eqref{eq:AllPolZConstraint} and the causality requirement are rewritten for compactness as a set of admissible policies $(D_i,e_i)$ parameterized by the current state $x^i_0$:
\begin{displaymath}
\mathcal{F}_i(x^i_0) = \left\{ \! (D_i, e_i)  \left| \begin{array}{l} [D_i]_{l,m} = 0,\, \hspace{2.2cm} \forall m > l  \vspace{0.2cm} \\   \left[ \!\!\!  \begin{array}{c} Ax^i_0\! +\! B_i(D_i\delta + e_i) \\ D_i\delta + e_i \\ \delta \end{array} \!\!\! \right] \! \in \! \mathcal{Z}_i,\forall \delta \in \Delta \! \end{array} \!\!\!\! \right. \right\}
\end{displaymath}
This leads to the following rewriting of problem \eqref{eq:OptOverAllPolicies} in terms of nominal schedules $\{e_i\}_{i=1}^{N_{\rm p}}$ and policy matrices $\{D_i\}_{i=1}^{N_{\rm p}}$:
\begin{subequations} \label{eq:AffinePolOpt}
\begin{IEEEeqnarray}{L}
\min_{(D_i,e_i) \in \mathcal{F}_i(x^i_0)} \quad  \sum_{i=1}^{N_{\rm p}} \tilde{J}_i(x^i_0,D_i,e_i) \label{eq:DistrObjFunction}\\
\text{s.t.}\,  \sum_{i=1}^{N_{\rm p}} r_i + G_i \delta + C_i(A_ix^i_0 + B_i(D_i\delta + e_i)) = 0 \, , \vspace{-0.2cm} \nonumber \\
\hspace{6.2cm} \forall \delta \in \Delta \, , \label{eq:DistrJointMatchConstraint}\\
\hspace{0.5cm} \sum_{i=1}^{N_{\rm p}} \Gamma_i (r_i + G_i \delta + C_i(A_ix^i_0 + B_i(D_i\delta + e_i))) \leq \overline{p} \, , \vspace{-0.2cm} \nonumber \\
\hspace{6.2cm} \forall \delta \in \Delta \, .  \label{eq:DistrJointLineConstraint}
\end{IEEEeqnarray}
\end{subequations}

The objective function has been redefined as $\tilde{J}_i(x^i_0,D_i,e_i) := \mathbb{E}[J_i(\mathbf{x}^i,\mathbf{u}^i)]$ to reflect its dependence on $x_0^i$, $D_i$, and $e_i$.

The assumption of a positive semidefinite quadratic form \eqref{eq:CostFunction} for $J_i(\mathbf{x}^i,\mathbf{u}^i)$ allows the expectation over $\delta$ to be minimized straightforwardly in comparison to the arbitrary case, since only the moments $\mathbb{E}[\delta]$ and $\mathbb{E}[\delta \delta']$ are needed. Substitution from \eqref{eq:CostFunction} gives the following representation of the objective, which is convex:
\begin{IEEEeqnarray*}{RL}
\tilde{J}_i(x^i_0,D_i,e_i) = & \mathbb{E}[J_i(A_i x^i_0 + B_i(D_i \delta + e_i),D_i \delta + e_i)] \\
= & f_i^x{}'(A_ix^i_0 + B_ie_i) + f_i^u{}'e_i + x^i_0{}'A_i'H_i^xB_ie_i \\
+ & \tfrac{1}{2}x^i_0{}'A_i'H_i^xA_ix^i_0 + \tfrac{1}{2}e_i'(B_i'H_i^xB_i + H_i^u)e_i\\
+ & (f_i^x{}'B_i + f_i^u{}' + x^i_0{}'A_i'H_i^xB_i \\
&\hspace{1.5cm} + e_i'B_i'H_i^xB_i + e_i'H_i^u)D_i\mathbb{E}[\delta] \\
+ & \tfrac{1}{2} \langle D_i'(B_i'H_i^xB_i + H_i^u)D_i,\mathbb{E}[\delta \delta'] \rangle + c_i 
\end{IEEEeqnarray*}

Since in many cases the reference predictions $r_i$ will be chosen with $\mathbb{E}[\delta]=0$, the corresponding term above generally cancels. Note that although the constant term $c_i$ makes no difference to the solutions $D_i$ and $e_i$, it is needed later in order to compare the costs of different approaches.

\subsection{Equivalent tractable reformulation}

Problem \eqref{eq:AffinePolOpt} cannot be solved directly because constraints \eqref{eq:DistrJointMatchConstraint} and \eqref{eq:DistrJointLineConstraint}, as well as the definition of $\mathcal{F}_i(x^i_0)$, apply for all $\delta \in \Delta$ and are therefore the intersection of an infinite number of constraints. To obtain a numerical solution they must be written in an equivalent finite form.

Since $\Delta$ is assumed to have an interior containing the origin, i.e.~no redundant components exist in $\delta$, it can be shown that constraint \eqref{eq:DistrJointMatchConstraint} is satisfied if and only if the following conditions hold:
\begin{subequations}
\begin{IEEEeqnarray}{RL}
\sum_{i=1}^{N_{\rm p}} (r_i + C_i A_i x^i_0 + C_i B_i e_i) & = 0\, ,\\
\smash{\sum_{i=1}^{N_{\rm p}}} (G_i + C_i B_i D_i) & = 0\, .\vspace{0.2cm}
\end{IEEEeqnarray}
\end{subequations}

Constraint \eqref{eq:DistrJointLineConstraint} and the sets $\mathcal{F}_i(x^i_0)$ can be written using a result due to Guslitser and others \cite{guslitser_uncertainty-immunized_2002, ben-tal_adjustable_2004}. Recalling that $\Delta = \{ \delta \, | \, S\delta \leq h \}$, the following equivalences hold for \eqref{eq:DistrJointLineConstraint}. An extra matrix variable $Z$ is introduced for the last equivalence, which uses strong duality in linear programming (see \cite{goulart_optimization_2006}, Example 7).
\allowdisplaybreaks
\begin{IEEEeqnarray*}{C}
\sum_{i=1}^{N_{\rm p}} \Gamma_i (r_i + G_i \delta + C_i(A_ix^i_0 + B_i(D_i\delta + e_i)))  \leq \overline{p}, \, \forall \delta \in \Delta \\
\Updownarrow  \\
\hspace{-2.2cm}\max_{\delta \in \Delta} \sum_{i=1}^{N_{\rm p}} \Gamma_i(G_i + C_i B_i D_i)\delta + \sum_{i=1}^{N_{\rm p}} \Gamma_i C_i B_i e_i \\
\hspace{4.2cm} \leq \overline{p} - \sum_{i=1}^{N_{\rm p}} \Gamma_i(r_i + C_i A_i x^i_0) \\
\Updownarrow  \\
\hspace{-0.1cm} \left\{ \!\!\! \begin{array}{l} \exists Z \! : \! Z'h \! + \! \sum_{i=1}^{N_{\rm p}} \Gamma_i C_i B_i e_i  \leq \overline{p} - \! \sum_{i=1}^{N_{\rm p}} \Gamma_i(r_i + C_i A_i x^i_0),  \vspace{0.2cm} \\
 \sum_{i=1}^{N_{\rm p}} \Gamma_i(G_i + C_i B_i D_i) \! = \! Z'S,  \, \text{and}\, Z \geq 0\,\, \text{element-wise.} \end{array} \!\!\! \right\}
\end{IEEEeqnarray*}

Similarly, it can be shown that the sets $\mathcal{F}_i(x^i_0)$ can be rewritten in finite form as follows, starting from definition \eqref{eq:ZSetDef} and introducing extra matrix variables $Y_i$ of appropriate dimension:
\begin{IEEEeqnarray*}{L}
 \mathcal{F}_i(x^i_0)\!  =  \! \left\{ \! (D_i, e_i) \left| \! \begin{array}{l} [D_i]_{l,m} = 0,\, \forall m > l \\ \exists Y_i \geq 0 \!: (T_iB_i + U_i)D_i + V_i = Y_i'S,\\  
T_iA_ix^i_0 + (T_iB_i + U_i)e_i + Y_i'h \leq w_i \end{array} \!\!\! \right. \right\}
\end{IEEEeqnarray*}




These changes lead to the following tractable representation of optimization \eqref{eq:AffinePolOpt}:
\begin{subequations} \label{eq:DistrOpt}
\begin{IEEEeqnarray}{L}
\min_{Z \geq 0,\,(D_i,e_i) \in \mathcal{F}_i(x^i_0)} \quad  \sum_{i=1}^{N_{\rm p}} \tilde{J}_i(x^i_0,D_i,e_i)  \\
\text{s.t.}\,\,  \sum_{i=1}^{N_{\rm p}} (r_i + C_i A_i x^i_0 + C_i B_i e_i)  = 0\, , \label{eq:DistrConstr1}\\
\hspace{0.5cm} \sum_{i=1}^{N_{\rm p}} (G_i + C_i B_i D_i)  = 0\, , \label{eq:DistrConstr2}\\
\hspace{0.5cm} Z'h + \sum_{i=1}^{N_{\rm p}} \Gamma_i C_i B_i e_i  \leq \overline{p} - \sum_{i=1}^{N_{\rm p}} \Gamma_i(r_i + C_i A_i x^i_0)\, , \hspace{0.5cm}\label{eq:DistrConstr3}\\
\hspace{0.5cm} \sum_{i=1}^{N_{\rm p}} \Gamma_i(G_i + C_i B_i D_i)  = Z'S\, . \label{eq:DistrConstr4}
\end{IEEEeqnarray}
\end{subequations}

In summary, \eqref{eq:DistrConstr1} states that the nominal schedules of power output changes $e_i$ should track the base prediction; \eqref{eq:DistrConstr2} states that the rules $D_i$ used by the participants should together track any error vector $\delta \in \Delta$; \eqref{eq:DistrConstr3} and \eqref{eq:DistrConstr4} ensure that line current limits should not be exceeded for any $\delta \in \Delta$.

After solving \eqref{eq:DistrOpt} the state and input trajectories $\mathbf{x}^i$ and $\mathbf{u}^i$ for a particular prediction error $\delta$ can be computed by substituting the solution $(D_i,e_i)$ back into \eqref{eq:AffinePolDef} and \eqref{eq:Dynamics}.

\subsection{Computational requirements}

Problem \eqref{eq:DistrOpt} is a quadratic program, which in principle can be solved even where many thousands of variables and constraints are present. We now comment on the size and structural properties of the problem. Each vector $e_i$ has $T$ elements, each matrix $D_i$ has $N_\delta T^2$ elements (neglecting the fact that some of these are constrained to be zero), and each matrix $Y_i$ (introduced by definition of the sets $\mathcal{F}_i (x^i_0 )$) has $ql_i$ elements, recalling that $q$ is the number of constraints defining $\Delta$ and $l_i$ is the number of constraints defining $\mathcal{Z}_i$. The matrix $Z$ has $2qLT$ elements. Therefore the total number of primal optimization variables is $N_{\rm p} (T + N_\delta T^2 + q \overline{l}) + 2qLT$, where $\overline{l} = \tfrac{1}{N_{\rm p}}\sum_{i=1}^{N_{\rm p}}l_i$. The main computational cost arises from matrices $D_i$ which together have $N_{\rm p} N_\delta T^2$ elements. The problem size therefore grows quadratically with the time horizon, and already reaches the thousands for modest parameter choices. In contrast, the number of decision variables needed to operate reserves in a way more comparable with existing mechanisms (see Section \ref{sec:CostComparison}) grows only linearly with the time horizon.

However, for computational purposes the structure of the problem is as important as the size. The form of \eqref{eq:DistrOpt} exhibits a convenient linear coupling between participants $i= 1,\ldots,N_{\rm p}$ in both the cost function and the constraint set, and in fact lends itself to solution via a large-scale distributed solution method, such as the recently-revived Alternating Direction Method of Multipliers \cite{boyd_distributed_2010}. For the examples reported in this paper, though, the resulting optimization problems were still manageable enough for a centralized solution. 

\section{Implementation of policy-based reserves} \label{sec:MarketSol}

Today, electrical reserves are often allocated through market means (see \cite{entso-e_survey_2012} for a survey) as a way of shifting the onus to provide them efficiently onto market participants. However, a parallel trend towards sophisticated centralized optimizations carried out by system operators has also arisen \cite{ott_evolution_2010}, in cases where real-time decisions such as unit commitment cannot be settled adequately by market means. Both market and non-market schemes could be envisaged for the provision of policy-based reserves. Regardless of whether a market is used, though, appropriate payments to participants for the reserve service must still be determined. This section derives these payments, then considers what form a traded reserve policy product may take, and how it may be executed.

\subsection{Efficient prices for policies} \label{sec:PriceDer}

We now show that efficient market prices exist for reserve policies and that these are an exact analogue of standard Locational Marginal Prices (LMPs), which arise in electricity markets from network congestion \cite{stoft_power_2002}. A partial Lagrangian of problem \eqref{eq:DistrOpt} (keeping the constraints $(D_i,e_i) \in \mathcal{F}_i(x^i_0)$ and $Z \geq 0$ but relaxing all others), $\mathcal{L}(Z,D_1,e_1,\ldots,D_{N_{\rm p}},e_{N_{\rm p}},\lambda,\Pi,\nu,\Psi)$, can be formed by introducing the multipliers $\lambda \in \mathbb{R}^T$ for violation of constraint \eqref{eq:DistrConstr1}, $\Pi \in \mathbb{R}^{T \times N_\delta T}$ for \eqref{eq:DistrConstr2}, $\nu \in \mathbb{R}_{\geq 0}^{2LT}$ for \eqref{eq:DistrConstr3}, and $\Psi \in \mathbb{R}^{2LT \times N_\delta}$ for \eqref{eq:DistrConstr4}. 
After some rearrangement, this partial Lagrangian can be rewritten in separable form as
\begin{IEEEeqnarray}{L}
\mathcal{L}(Z,D_1,e_1,\ldots,D_{N_{\rm p}},e_{N_{\rm p}},\lambda,\Pi,\nu,\Psi) = \nonumber \\
\hspace{2.5cm} \sum_{i=1}^{N_{\rm p}} (\tilde{J}_i(x^i_0,D_i,e_i) - \lambda_i'e_i - \langle \Pi_i,D_i \rangle) \nonumber \\
                       \hspace{2.5cm} +  \nu'Z'h - \langle \Psi, Z'S \rangle + f(\lambda,\Pi,\nu,\Psi)  \hspace{1cm}
\end{IEEEeqnarray}
where $f$ is constant with respect to the primal variables $D_i$, $e_i$, and $Z$. The Lagrange multipliers from optimization problems are commonly interpreted as prices, and here we have in effect defined two prices. The first is a nodal power price \[ \lambda_i :=  -B_i'C_i'(\lambda + \Gamma_i'\nu) \, ,\] consisting of a global component depending on $\lambda$ and a local component (induced by any line congestion present) depending on $\nu$. This agrees with the standard derivation of locational marginal prices (LMPs) in optimal power flow theory. The second is a matrix of reserve policy prices \[ \Pi_i :=  -B_i' C_i' (\Pi + \Gamma_i'\Psi)\, ,\] whose entries are the marginal value of each element of $D_i$. Note that the elements of $\Pi_i$ above the main diagonal are not required, since these are related to entries that would anyway result in non-causal responses to prediction errors, were they to be different from zero.

It can be shown that $C_iB_i = I$ in most cases (see the Appendix), so that usually $\lambda_i = -(\lambda + \Gamma_i'\nu)$ and $\Pi_i = -(\Pi + \Gamma_i'\Psi)$. The minus signs are a result of the sign convention used to write constraints \eqref{eq:DistrConstr1} and \eqref{eq:DistrConstr2}. The identical form of $\lambda_i$ and $\Pi_i$ suggests that optimal prices for reserve policies exhibit the same locational variation as LMPs.

Because the Lagrangian is separable, the terms $\lambda_i'e_i$ and $\langle \Pi_i,D_i \rangle$ can be identified as the efficient payments to each participant $i$ for committing to nominal plan $e_i$ and reserve policy $D_i$. The payments represent the money transfers that would result from a market mechanism that solved problem \eqref{eq:DistrOpt} efficiently. Under such a scheme, the expected profit made by participant $i$ would be equal to $-\tilde{J}_i(x^i_0,D_i,e_i) + \lambda_i'e_i + \langle \Pi_i,D_i \rangle$.


\subsection{Policy-based reserve products} \label{sec:ReserveProducts}

A natural question is what physical commitment a participant makes when implementing a given policy $(D_i,e_i)$. Clearly the nominal part $e_i$ is the schedule participant $i$ will follow if predictions turn out to have been made with perfect accuracy, i.e.~$\delta=0$. The interpretation of $D_i$ is more subtle and can be termed in two ways, which will be described here. In this subsection the analysis is given for a single uncertainty source, i.e. $N_\delta = 1$; the case for $N_\delta > 1$ is analogous. 

Firstly, each \emph{row} $l$ of $D_i$ can be read as the rule the participant must follow to construct its input $u_l^i$, which for the realized errors $\delta_0,\ldots,\delta_l$ determines the power it injects at time $l+1$, so that $u_l^i = [e_i]_l + \sum_{m=0}^l [D_i]_{l,m} \delta_m$.

Secondly, each \emph{column} $m$ of $D_i$ (reading down from the element on the diagonal) is, in control terminology, the planned \emph{impulse response} $g(\delta_{m+1})$ of participant $i$ to a unit prediction error $m$ steps from the current time (recalling that $\delta_{m+1}$ is revealed just before $u^i_m$ is applied), so that $g(\delta_{m+1}) = [  [D_i]_{m,m}, [D_i]_{m+1,m}, \ldots, [D_i]_{T-1,m} ]$.

\subsubsection*{Example product} Consider the example of the decision rule governing $u^i_l$, which sets the power participant $i$ supplies during step $l+1$ of the planning horizon. For given choices of $[D_i]_{l,0}$, $[D_i]_{l,1}$, up to $[D_i]_{l,l}$ agreed, the payment (using the price notation from Section \ref{sec:PriceDer}) made to participant $i$ for the reserve service would be $\sum_{m=0}^{l}[\Pi_i]_{l,m}[D_i]_{l,m}$ and the payment for scheduled power $[e_i]_l$ would be  $[\lambda_i]_l[e_i]_l$. An analogous product could be sold based on the column-wise reading of $D_i$, in which case the participant would be selling an \emph{a posteriori} response to errors.

Existing reserve mechanisms could be modelled by matrices $D_i$ for which only the main diagonal is populated. In secondary reserves provided by conventional generators, feedback controllers use the frequency deviation, in the form of the Area Control Error (ACE) signal, to adjust their power outputs to follow load mismatches. For a unit deviation from the net load reference, each generator will end up with an offset that depends on its controller parameters. For such a unit load deviation at time $l$, this offset is in our notation exactly the matrix entry $[D_i]_{l,l}$.

Therefore, a simple way of comparing policy-based reserves to existing mechanisms is to restrict the structure of the matrices $D_i$ accordingly and inspect the results. This is demonstrated in Section \ref{sec:NumEx}.

\subsection{Real time operation} \label{sec:RTOperation}

For continuous operation of the power system, it is necessary to choose new policies periodically, because the policies only apply for the following $T$ steps at the time when they are chosen. Furthermore the optimal policies are a function of the current state of the system, and it may be attractive to choose new policies early in the light of new forecast information. Therefore a systematic way of choosing policies repeatedly is required. Two possible schemes for such ``closed-loop'' operation are:
\begin{itemize}
	\item \textbf{Batchwise}: Set policies for a horizon of $T$ steps, let all $T$ steps play out, then choose a new batch of policies once the current ones have expired.
	\item \textbf{Receding horizon}: Let only the first step play out, then immediately update the policies for the next $T$ steps.
\end{itemize}
Clearly, a middle road between these two options also exists, namely letting some of the $T$ steps play out, then discarding the remainder and choosing new policies.

The receding horizon scheme presents another choice --- whether to honour previous policies in the new choices for $D_i$ (respectively $e_i$). This would be done by shifting the matrices (resp.~vectors) upward and leftward (resp.~upward) by one element, then optimizing over the bottom row of elements (resp.~element) to obtain new policies. These new policies would be feasible with respect to $\mathcal{F}_i(x_0^i)$ as well as the global constraints \eqref{eq:DistrConstr1}-\eqref{eq:DistrConstr4} as long as $\Delta$ is defined such that the bounds on the uncertainty grow monotically along the prediction horizon (for conciseness this is stated without proof here).

An alternative to this is to reject the previous policies and choose fresh values for $D_i$ and $e_i$ at every step. This option is attractive because it allows new (presumably better) policies to be chosen in the light of new information as soon as it becomes available. However its drawback is that due to the repeated renewal process, only the $(0,0)$ block of the matrix $D_i$ ever gets used, and no policies beyond the first row ever appear to be implemented. This conflicts with the ideas developed in Section \ref{sec:ReserveProducts}, in that contracts for policy-based reserves would be agreed but then never called on. It is important to note, however, that the optimal values of $[e_i]_0$ and $[D_i]_{0,0}$ would be affected strongly by the costs and constraints modelled for steps $1$ to $T-1$, even though those later steps would never be realized. This means that under this scheme, any cost savings from employing reserve policies are ultimately to be found in the more intelligent choice of $[e_i]_0$ and $[D_i]_{0,0}$. This effect is reported in Section \ref{sec:NumEx}.

Rejecting previous policies at each time step makes choosing correct payments to market participants difficult. One simple way of overcoming this would be instead to define the new policies as deviations from shifted versions of the existing previous policies, with payments settled additively. Under such a scheme there would be $T$ payments for the reserve action undertaken at a given time, resulting from the superposition of policy choices made over the last $T$ steps.

In the presence of a large renewable infeed, bounds on the prediction errors may change significantly over the course of a prediction horizon; by the time the prediction horizon has nearly been played out, prediction error bounds for the last few steps are likely to be much smaller than those assumed at the start of the horizon, when the policies were chosen. Therefore it seems logical to prefer a rolling system of adjustments to policies, of the kind described above, rather than a batchwise approach.


\section{Numerical case study} \label{sec:NumEx}

We applied policy-based reserves to a standard test network that was adapted to include a large share of wind power. Policies were recomputed in a receding horizon fashion (the second method described in Section \ref{sec:RTOperation}) over a three-day simulated period in order to assess the reserve costs incurred. The test was repeated 50 times with different realizations of the random wind infeed. The effects of restrictions on the structure of the matrices $D_i$ on the cost of reserves were measured, leading to observations on the cost savings facilitated by recourse. The optimization problems were solved using CPLEX \cite{cplex_11.0_2007}.

The network used was a modification of the 39 bus network described in Appendix A of \cite{pai_energy_1989}, and shown in Fig.~\ref{fig:39BusNetwork}. This network contains 7 thermal generators, 2 storage units, 19 loads, and 3 wind farms (which replace 3 of the original 10 generators). We assumed that the generators represent the plants that had been selected for use via an earlier unit commitment decision. Generators have fuel costs represented in $f_i^u$ and $H_i^u$, and ramping costs represented in $H_i^x$. Storage units have a penalty for deviating from their midpoint, represented in $c_i$, $f_i^x$, and $H_i^x$. The parameters in terms of the definitions in Section \ref{sec:PartModel} are described in Table \ref{tab:Participants}. 

The daily pattern of load variation shown with the thick black line in Fig.~\ref{fig:PowerTraces} was taken from data for total UK national electrical consumption\footnote{available at \texttt{http://www.bmreports.com}} on 14${}^\text{th}$ September 2012, normalized to the size of each load modelled. Peak load was 6.097 GW, and wind power provided 3 GW at maximum, with an expected energy share of 29.0\%. Line flows from bus 16 to 15 and from bus 16 to 17 were restricted to 1000 MW. The load sizes $p_\text{nom}$ in Table \ref{tab:Participants} are the nominal values described in \cite{pai_energy_1989}. The three-day simulation period was divided into 288 fifteen-minute steps ($\tau = 0.25$ hrs), and the horizon length was $T=8$.

\begin{figure}[t]
  \centering
  \includegraphics[width=8cm]{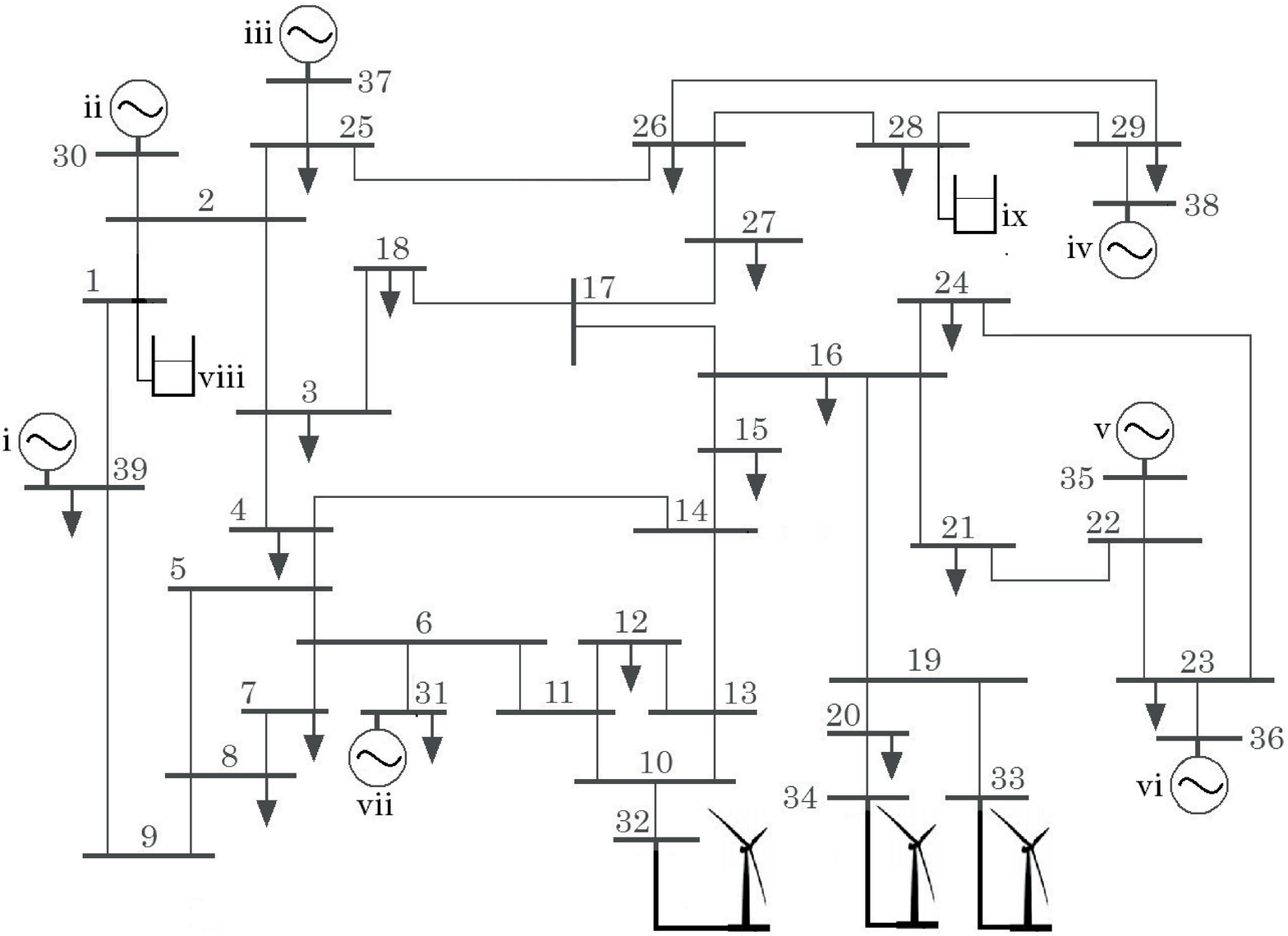} 
	\caption{39 bus test network from \cite{pai_energy_1989}, with wind infeed replacing thermal generators at nodes 32, 33, and 34, and with added storage units at nodes 1 and 28.}
	\label{fig:39BusNetwork}

\end{figure}

\begin{table}[t]
	\centering
		\caption{Parameters of elastic and inelastic participants}
		\begin{tabular}{|c|c|c|c|c|c|}
		\hline
		\multicolumn{6}{|l|}{\textbf{Thermal generators}, $i=1,\ldots,7$: \hspace{3cm} } \\
		\multicolumn{6}{|l|}{States: $\left[ \!\!\! \begin{array}{c} \text{Current output (MW)} \\ \text{Previous output (MW)} \end{array} \!\!\! \right]$,  $\tilde{A}_i = \left[ \!\!\! \begin{array}{cc} 0 & 0 \\ 1 & 0 \end{array} \!\!\! \right]$,  $\tilde{B}_i = \left[ \!\!\! \begin{array}{c} 1 \\ 0  \end{array} \!\!\! \right]$,  } \\
		\multicolumn{6}{|l|}{$\tilde{C}_i = \left[ 1 \, 0  \right]$, $\tilde{f}_i^x = \left[ \begin{array}{c} 0 \\ 0 \end{array} \right]$,  $\tilde{H}_i^x = \left[ \begin{array}{cc} \alpha & -\alpha \\ -\alpha & \alpha \end{array} \right]$, $\tilde{c}=0$, $x^i_0 = \left[ \!\!\! \begin{array}{c} p_0 \\ p_0  \end{array} \!\!\! \right]$ } \\
		\multicolumn{6}{|l|}{Constraints, $\forall \, k$: $0 \leq [x_k^i]_1 \leq p_\text{max}$, $0 \leq [x_k^i]_2 \leq p_\text{max}$,}\\
		\multicolumn{6}{|r|}{$0 \leq u_k^i \leq p_\text{max}$ } \\
		\hline
		$i$ & $\qquad \tilde{f}_i^u \qquad$ & $ \qquad \tilde{H}_i^u \qquad $ & $\qquad \alpha \qquad $ & $\quad p_\text{max}\quad $ & $p_0$\\
		\hline
		1 & 20 & 0.020 & 1.0 & 1800 & 400 \\
		2 & 20 & 0.080 & 0.1 & 450  & 100 \\
		3 & 20 & 0.037 & 0.1 & 972  & 216 \\
		4 & 20 & 0.024 & 1.0 & 1494 & 332 \\
		5 & 20 & 0.031 & 0.1 & 1170  & 260 \\
		6 & 20 & 0.036 & 0.1 & 1008  & 224 \\
		7 & 20 & 0.200 & 0.1 & 180  & 40 \\
	  \hline
	  \end{tabular}
		
	  \vspace{0.1cm}
	  \begin{tabular}{|c|c|c|c|c|}
	  \hline
		\multicolumn{5}{|l|}{\textbf{Storage units}, $i=8,9$:} \\
		\multicolumn{5}{|l|}{States: $\left[ \!\!\! \begin{array}{c} \text{Current output (MW)} \\ \text{Previous output (MW)} \\ \text{Storage level (MWh)}\end{array} \!\!\! \right]$,  $\tilde{A}_i = \left[ \!\!\! \begin{array}{ccc} 0 & 0 & 0 \\ 1 & 0 & 0 \\ 0 & 0 & 1 \end{array} \!\!\! \right]$,  $\tilde{B}_i = \left[ \!\!\! \begin{array}{c} 1 \\ 0 \\ -\tau \end{array} \!\!\! \right]$,  } \\
		\multicolumn{5}{|l|}{$\tilde{C}_i = \left[ 1 \, 0 \, 0 \right]$, $\tilde{f}_i^x = \left[ \begin{array}{c} 0 \\ 0 \\ -2 \gamma \frac{s_\text{max}}{2} \end{array} \right]$,  $\tilde{H}_i^x = \left[ \begin{array}{ccc} 0 & 0 & 0 \\ 0 & 0 & 0 \\ 0 & 0 & 2\gamma \end{array} \right]$,     } \\
		\multicolumn{5}{|l|}{$f^u_i = 0$, $\tilde{H}^u_i = 0$, $x^i_0 = \left[ 0 \,\, 0 \,\, s_0 \right]'$, $\tilde{c} = \gamma(\frac{s_\text{max}}{2})^2$} \\
		\multicolumn{5}{|l|}{Constraints, $\forall \, k$: $-p_\text{max} \leq [x_k^i]_1 \leq p_\text{max}$, $-p_\text{max} \leq [x_k^i]_2 \leq p_\text{max}$,}\\
		\multicolumn{5}{|r|}{ $0 \leq [x_k^i]_3 \leq s_\text{max}$, $-p_\text{max} \leq u_k^i \leq p_\text{max}$ } \\
		\hline
		$i$ & $\qquad s_\text{max} \qquad$ & $\qquad \gamma \qquad$ & $\qquad p_\text{max} \qquad$ & $s_0$ \\
		\hline
		8 & 1000 & 0.01 & 200 & 500 \\
		9 & 1000 & 0.01 & 200 & 500 \\
	  \hline
	  \end{tabular}
		
	  \vspace{0.1cm}
	  \begin{tabular}{|r|}
	  		\hline		
	  		\textbf{Wind farms}	[\emph{Node:} $\tilde{G}_i$] \hspace{1.5cm} \emph{32:} $[2\,\, 0]$, \emph{33:} $[1 \,\, 1]$, \emph{34:} $[0 \, \, 2]$ \\
		\hline
		 \textbf{Loads} [\emph{Node:} $p_\text{nom}$] \hspace{0.6cm} \emph{3:} 322, \emph{4:} 500, \emph{7:} 233.8, \emph{8:} 522, \emph{12:} 7.5,\\
		\emph{15:} 320, \emph{16:} 329, \emph{18:} 158, \emph{20:} 628, \emph{21:} 274, \emph{23:} 247.5, \emph{24:} 308.6, \\
	  \emph{25:} 224, \emph{26:} 139, \emph{27:} 281, \emph{28:} 206, \emph{29:} 283.5, \emph{31:} 9.2, \emph{39:} 1104 \\
		\hline
	  \end{tabular}
	\label{tab:Participants}

\end{table}

\subsection{Uncertainty model}

Uncertainties in the system were assumed to originate only in the random wind power availability (loads were assumed to be predicted exactly, though our method could equally be used to account for load uncertainty). The wind farm output was driven by the following first order random process model with saturation: 
\begin{equation} \label{eq:RandomProcess}
 q_{k+1} = \min \{ \max \{q_\text{min}, q_k + \beta_k \}, q_\text{max} \} \, , 
\end{equation}
where $q_k \in \mathbb{R}^{N_\delta}$ denotes the state of the uncertainty model, and $\beta_k$ is sampled at each step $k$ from a multivariate normal distribution with variance $\Sigma \in \mathbb{R}^{N_\delta \times N_\delta}$, truncated to bounds $A_\beta \beta_k \leq b_\beta$. Note that more elaborate wind models exist \cite{papaefthymiou_mcmc_2008}, but we have used a simpler model here for demonstration purposes.

Defining $q := [q_1' \ldots q_T']'$ as the random future evolution of $q$ from current state $q_0$, the nominal predictions of wind farm power output $r_i$ were mapped linearly from $\mathbb{E}[q]$, so that $r_i = G_i \mathbb{E}[q]$, where $G_i$ is the same matrix as that described in Section \ref{sec:PartModel}. The prediction error was defined as $\delta := q - \mathbb{E}[q]$ so that $\mathbb{E}[\delta] = 0$. The prediction error covariance is then $\mathbb{E}[\delta \delta'] = \mathbb{E}[(q - \mathbb{E}[q])(q - \mathbb{E}[q])']$. This was supplied together with the references $r_i$ and the current system state $x^i_0$ as inputs to optimization problem \eqref{eq:DistrOpt}. 

The uncertainty set $\Delta$ was recomputed as a function of the current state, to reflect the fact that prediction errors are bounded relative to the nominal predictions by the wind farm power output limits and by the bounds assumed on $\beta_k$, the change in wind power availability at each step. Its new parameters $S$ and $h$ were then supplied to \eqref{eq:DistrOpt}. In addition, at every simulation time step, an aggregation of 20,000 Monte Carlo runs was used to produce $T$-step nominal predictions $r_i$ for each wind farm, as well as estimates of $\mathbb{E}[\delta \delta']$. 

\begin{figure}[t]
  \centering
  \includegraphics[width=13.6cm]{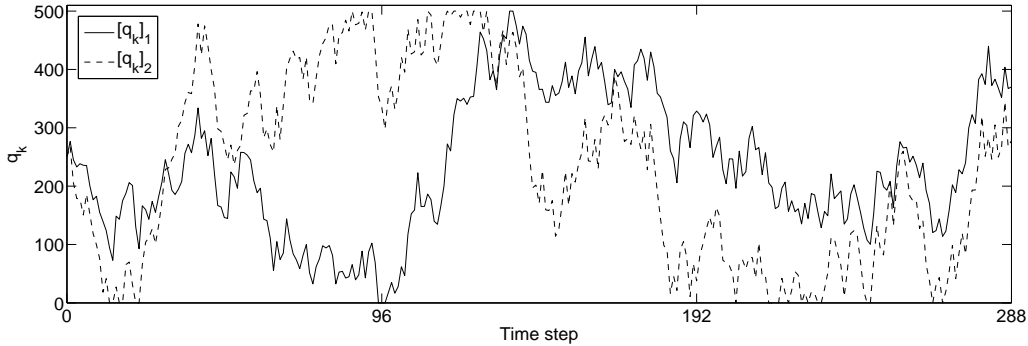}
	\caption{Example output of correlated uncertainty model $q_k$, which drives the wind farm outputs in the case study. The solid line represents $[q_k]_1$ and the dotted line $[q_k]_2$.}
	\label{fig:WindTrace}

\end{figure}

The three wind farms in this case study are driven by two temporally and spatially correlated sources of uncertainty, with parameters $\Sigma = [{}^{2000}_{2400} \, {}^{2400}_{2000}]$, $q_\text{min} = [{}^0_0 ]$, $q_\text{max} = [{}^{500}_{500} ]$, $A_\beta = [{}^{\, I}_{-I} ]$, $b_\beta = 80 \cdot \mathbf{1}$. The matrices $G_i$ for the wind farms are of the form $G_i = I_T \otimes \tilde{G}_i$, and the matrices $\tilde{G}_i$ are given in Table \ref{tab:Participants}. The output of the wind farm at node 34 is the mean of the two at nodes 32 and 33. The purpose of this formulation is to show that random power flows can be driven by a process with dimension lower than the number of nodes affected, either to save computational effort, or because the uncertainty is anyway difficult to model. An example time series $q_k$ generated in this way is shown in Fig.~\ref{fig:WindTrace}.

\subsection{Comparison of reserve costs} \label{sec:CostComparison}

For each of the 50 wind realization tests, three parallel models of the power system were driven by the same realizations of the random wind model described above. The observed operating costs under receding horizon control over the simulation period were then compared under the three schemes, which were defined as follows:

\begin{enumerate}
	\item Prescient case: Disturbances are known at the time the finite-horizon optimization is carried out. Matrices $D_i$ are therefore not needed, and nominal schedules $e_i$ track the power reference perfectly. This scheme, which results in the best attainable receding-horizon cost, is used as a point of comparison for the other two.
	\item Flexible-rate reserves: $[D_i]_{l,m} = 0$ for  $l \neq m$, for all $i$. This represents the best possible response to uncertainty without time coupling, and the optimization is over the elements of $e_i$ and the diagonal parts of $D_i$.
	\item Policy-based reserves: $[D_i]_{l,m} = 0$ for  $l < m$, for all $i$. This allows full use of the extra information that will be available at each time step when the reserve is deployed.
\end{enumerate}

%

The total operation cost was measured from the state and input values $x^i(t)$ and $u^i(t)$ realized by the elastic participants at each time step $t$ over the simulation period, 
\begin{align*}
&\sum_{t=1}^{288} \sum_{i=1}^{N_{\rm p}} \left[ \tilde{f}_i^x{}'x^i(t) + \frac{1}{2}x^i{}'(t)\tilde{H}_i^x x^i(t) + \tilde{f}_i^u{}'u^i(t-1) \,+ \right. \\
& \hspace{4cm}  \left. \frac{1}{2}u^i{}'(t-1)\tilde{H}_i^u u^i(t-1) + \tilde{c}_i \right]\, .
\end{align*} 
The inputs are indexed by $(t-1)$ because state $x(0)$ is given whereas input $u(0)$ must be chosen and determines $x(1)$, and so on.

\begin{figure}[t]

  \centering
  \includegraphics[width=13.6cm]{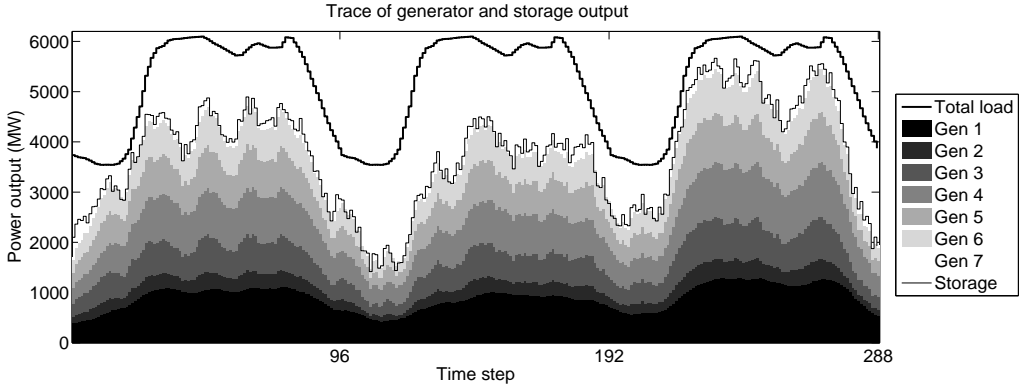} 

	\caption{Example of power output traces for one wind realization under Scheme 3. Generator outputs are plotted in stacked form, and the total power output of the two storage units is plotted using the thin black line relative to the top of the stack. The wind power injection is the difference between the total load (bold black line) and the sum of storage and generator infeeds (thin black line).}
	\label{fig:PowerTraces}

\end{figure}

An example of the power output traces for the generators is given in Fig.~\ref{fig:PowerTraces}. The cost results are shown in Table \ref{tab:StructureResults}. A cost of reserves is defined for Schemes 2 and 3 as the operation cost experienced minus the prescient cost (from Scheme 1). This represents the cost incurred in order to accommodate the uncertain wind infeed. Across the 50 runs, the cost of reserves decreased by an average of 38.4\% for full policies (Scheme 3) with respect to the best possible non-recourse reserve scheduler (Scheme 2). Average costs under Scheme 1 were $4.474 \times 10^7$ for the three-day test, and reserves were found to incur up to 1.10\% of additional total power generation costs in the tests under Scheme 2, and up to 0.73\% under Scheme 3. Computation times are also reported in Table \ref{tab:StructureResults}, and represent the average time needed to build and solve each finite-horizon optimization problem on an Intel Xeon E5-2670 2.60 GHz CPU. Note that Scheme 1 is far faster than the others because only open-loop schedules $e_i$ need to be computed.

Another series of tests was carried out allowing only the diagonal and immediate subdiagonal entries in the $D$-matrices to differ from zero. This results in a two-step policy, requiring fewer non-zero optimization variables than Scheme 3. The average savings in the cost of reserves, measured in the same way as for Scheme 3, were 32.4\%, indicating that on average most of the 38.4\% savings can be obtained from using just a two-step policy rather than a full policy (in this case 8 steps).

\begin{table}[t]
	\centering
		\caption{Cost comparison for different structural restrictions on $D_i$}

		\begin{tabular}{|c|c|c|c|}
		\hline
		Scheme & Comp. & Average cost increase & Avg.~reserve cost \\
		& time & over Scheme 1 & vs.~Scheme 2 \\
		\hline
		1. Prescient	   & 34 ms      & ---      & 0 \% \\
		2. Diagonal      & 3843 ms      & 0.825 \% & := 100.0 \% \\
		3. Full policy   & 4424 ms      & 0.511 \% & 61.6 \% \\
		\hline
		\end{tabular}
	\label{tab:StructureResults}

\end{table}

\subsection{Sensitivity analysis}

\begin{figure}[t]
  \centering
  \includegraphics[width=13.6cm]{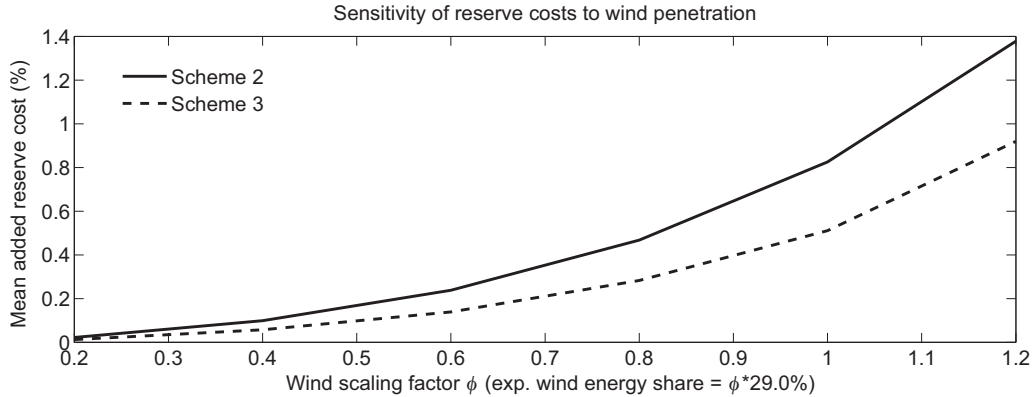} 
	\caption{Added reserve cost percentages under different wind penetration ratios $\phi$, averaged over 50 three-day simulation runs.}
	\label{fig:Sensitivity}
\end{figure}

The question arises whether the cost savings reported for affine reserve policies depend on the quantity of wind energy present. To test this the case study was repeated under wind realizations driven by different-sized instances of the random process \eqref{eq:RandomProcess}, scaling $q_\text{max}$ and $b_\beta$ by a factor $\phi \in [0.2,1.2]$, and $\Sigma$ by $\phi^2$ since $\Sigma$ represents the variance of a linearly-scaled quantity. The results, for 50 runs each, are shown in Table \ref{tab:Sensitivity} and plotted in Fig.~\ref{fig:Sensitivity} (note that for $\phi > 1.2$ infeasibility arises from the fact that wind power is not curtailed). The expected proportion of load energy supplied by wind is $\phi \cdot 29.0\%$.

The additional percentage cost of accommodating the uncertainty increased sharply (apparently more than quadratically) as the wind capacity was increased. Using time-coupled policies, around 40\% of this could be offset as the wind share grew. Although the percentage saving decreased slightly for increasing $\phi$, the absolute saving still grew quadratically. Note that the large saving reported for $\phi = 0.2$ is measured relative to a tiny number.

\begin{table}[t]
	\centering
		\caption{Variation of results with installed wind power capacity}

		\begin{tabular}{|c|c|c|c|c|}
		\hline
		Scaling & Wind     & \multicolumn{2}{c|}{Average reserve costs} & Average reduction  \\
		$\phi$  & capacity & Scheme 2 & Scheme 3     & under Scheme 3 \\
		\hline
		0.2 &  0.6 GW & 0.022 \% & 0.012 \% & 64.0 \% \\
		0.4 &  1.2 GW & 0.099 \% & 0.057 \% & 43.2 \% \\
		0.6 &  1.8 GW & 0.238 \% & 0.139 \% & 42.3 \% \\
		0.8 &  2.4 GW & 0.468 \% & 0.283 \% & 39.9 \% \\
		1.0 &  3.0 GW & 0.825 \% & 0.511 \% & 38.4 \% \\
		1.2 &  3.6 GW & 1.378 \% & 0.919 \% & 33.7 \% \\
		\hline
		\end{tabular}
	\label{tab:Sensitivity}

\end{table}

\section{Conclusions} \label{sec:Conclusion}

This paper introduced the idea of policy-based reserves, motivated by the need to incorporate intermittent renewable energy sources into power systems intelligently. The approach uses robust optimization to find efficient, \emph{time-coupled} responses to errors in the prediction of uncertain load or supply, over a finite time horizon. Interpretations were developed for policy-based reserves as products that could be bought or sold on power markets, with associated prices. A case study demonstrated the operating principles and benefits of the approach.

The potential cost savings from reserve policies are determined by the time coupling of participant costs, constraints, and prediction errors. This is evident from the term $\tfrac{1}{2} \langle D_i'(B_i'H_i^xB_i + H_i^u)D_i,\mathbb{E}[\delta \delta'] \rangle$ in the cost function. Prediction error correlation in space and time results in non-zero off-diagonal entries in $\mathbb{E}[\delta \delta']$, which will change the contributions of individual entries of the policy matrices $D_i$ to the expected costs. It is not surprising, then, given the strong correlation of prediction errors found in the numerical case study, and the presence of time-coupled generation costs, that for a given finite-horizon optimization the costs were reduced when reserve policies were enabled. 

In receding horizon operation, new power schedules are found at every time step. This can be viewed as a form of recourse that is not taken into account for each finite horizon optimization. This leads to the intuition that some of the apparent value of optimizing over full lower-triangular policies (rather than diagonal ones) may disappear once the optimization is repeated in receding horizon fashion. However our results show that cost savings indeed remain. It was also observed that much of the cost reduction can be gained without using all subdiagonal entries in the matrices $D_i$.

The interpretation of the reserve cost defined in Section \ref{sec:CostComparison} is important. Scheme 1 represents a lower bound on the cost achievable under any possible receding horizon scheme, but it may in fact be impossible to find causal policies $\pi_i(\delta)$ offering the required robustness to prediction errors with costs approaching this value. Therefore the savings we report are in fact lower bounds on the true available savings. Recent theoretical performance bounds \cite{hadjiyiannis_efficient_2011, van_parys_infinite-horizon_2012} suggest that the optimality of our reserve policies would depend on the current power system state, but that this relationship is far from straightforward.

We employed a standard set of power system simplifications in this paper in order to allow tractable optimization problems to be formulated. For a real implementation, additional ex-post checks (possibly with some iteration), would be needed to confirm feasibility of the flows on the AC transmission grid, in a fashion similar to current system operator practice \cite{oneill_recent_2011}.


\bibliographystyle{IEEEtran} 
\bibliography{ReserveLDRs,momoh,bertsimas,ferc}

\appendix

\subsection*{Cases where $C_iB_i = 0$ for elastic participants:}

Recall that from the definitions for $C_i$ and $B_i$ given in Section \ref{sec:PartModel},
\begin{align*} C_iB_i & = \left[ I_T \otimes [1 \, 0_{1 \times (n_i - 1)}] \right] \left[ \!\! \begin{smallmatrix} \tilde{B}_i & 0 & \cdots & 0 \\ \tilde{A}_i\tilde{B}_i & \tilde{B}_i & \ddots & \vdots \\ \vdots & \ddots & \ddots & 0 \\ \tilde{A}_i^{T-1}\tilde{B}_i & \cdots & \tilde{A}_i\tilde{B}_i & \tilde{B}_i \end{smallmatrix} \!\! \right] \\
& = \left[ \begin{smallmatrix} [\tilde{B}_i]_1 & 0 & \cdots & 0 \\ 
\text{$[\tilde{A}_i\tilde{B}_i ]_1$} & [\tilde{B}_i]_1 & \ddots & \vdots \\ 
\vdots & \ddots & \ddots & 0 \\ 
\text{$[\tilde{A}_i^{T-1}\tilde{B}_i ]_1$} & \cdots & \text{$[\tilde{A}_i\tilde{B}_i ]_1$} & [\tilde{B}_i]_1 \end{smallmatrix} \right] \, .
\end{align*}
As described in Section \ref{sec:PartModel}, we assume that the first state (the participant's net power infeed) is set directly by the input decision at the previous time step, and not as a function of current states. Therefore the first row of $\tilde{A}_i$ contains only zeroes, and $[\tilde{B}_i]_1 = 1$. From this it can easily be shown that $C_iB_i = I_T$, from which the property follows immediately.
\end{document}